\documentclass[reqno]{amsart}

\usepackage{amsmath}
\usepackage{amssymb}
\usepackage{amsthm}
\usepackage{xspace}
\usepackage{enumerate}

\theoremstyle{plain}
\newtheorem*{thm}{Theorem}
\newtheorem{lem}{Lemma}

\theoremstyle{remark}

\let\myskip=\medskip

\def\definebb#1=#2.{\def#1{{{\mathbb #2}^{\vphantom{x}}}}}
\definebb \c=C.
\definebb \r=R.
\definebb \z=Z.
\definebb \q=Q.

\def\chp{H_{\c}^2}
\def\rhp{H_{\r}^2}

\def\ppp{$(p_1,p_2,p_3)$}
\def\Gappp{\Ga(p_1,p_2,p_3)}

\def\al{{\alpha}}
\def\Ga{{\Gamma}}
\def\ga{{\gamma}}
\def\si{{\sigma}}
\def\om{{\omega}}

\def\st{\,\,\big|\,\,}
\def\<{\langle}
\def\>{\rangle}
\def\ie{i.e.\xspace}

\let\ge=\geqslant
\let\le=\leqslant

\def\defit{\it}

\DeclareMathOperator{\PU}{PU}
\DeclareMathOperator{\SU}{SU}
\let\Re=\undefined \DeclareMathOperator{\Re}{Re}
\let\mod=\undefined \DeclareMathOperator{\mod}{mod}

\begin{document}

\author[Anna Pratoussevitch]{Anna Pratoussevitch}
\address{Department of Mathematical Sciences\\ University of Liverpool\\ Peach Street \\ Liverpool L69~7ZL}
\email{annap@liv.ac.uk}

\title{Non-Discrete Complex Hyperbolic Triangle Groups of Type $(m,m,\infty)$}

\begin{date}  {\today} \end{date}

\thanks{Research partially supported by NSF grant DMS-0072607 
and by SFB 611 of the DFG}

\begin{abstract} 
In this note we prove that a complex hyperbolic triangle group of type $(m,m,\infty)$,
\ie a group of isometries of the complex hyperbolic plane,
generated by complex reflections in three complex geodesics meeting at angles $\pi/m$, $\pi/m$ and~$0$,
is not discrete if the product of the three generators is regular elliptic.

 \end{abstract}










\maketitle


\section{Introduction}

We study representations of real hyperbolic triangle groups,
\ie groups generated by reflections in the sides of triangles in~$\rhp$,
in the holomorphic isometry group $\PU(2,1)$ of the complex hyperbolic plane~$\chp$.

\myskip
For the basic notions of complex hyperbolic geometry, especially for the complex hyperbolic plane~$\chp$,
see for example section~2 in~\cite{P}.
The general references on complex hyperbolic geometry are~\cite{G99,P03}.

\myskip
We use the following terminology:
A {\defit complex hyperbolic triangle\/} is a triple $(C_1,C_2,C_3)$ of complex geodesics in~$\chp$.
If the complex geodesics $C_{k-1}$ and $C_{k+1}$ meet at the angle $\pi/p_k$
we call the triangle $(C_1,C_2,C_3)$ a {\defit \ppp-triangle}.

\myskip
We call a subgroup of $\PU(2,1)$ generated by complex reflections $\iota_k$ in the sides~$C_k$
of a complex hyperbolic \ppp-triangle $(C_1,C_2,C_3)$ a {\defit \ppp-triangle group\/}.
A {\defit \ppp-representation\/} is a representation of the group
$$\Gappp=\<\ga_1,\ga_2,\ga_3\st\ga_k^2=(\ga_{k-1}\ga_{k+1})^{p_k}=1~\hbox{for all}~k\in\{1,2,3\}\>,$$
where $\ga_{k+3}=\ga_k$, and the relation $(\ga_{k-1}\ga_{k+1})^{p_k}=1$ is to omit for $p_k=\infty$,
into the group $\PU(2,1)$,
given by taking the generators~$\ga_k$ of~$\Gappp$ to the generators~$\iota_k$ of a \ppp-triangle group.

\myskip
We prove in this paper the following result:

\begin{thm}
An $(m,m,\infty)$-triangle group is not discrete if the product of the three generators is regular elliptic.
\end{thm}

\myskip
More about the recent developments in the study of non-discrete complex hyperbolic triangle groups of type~$(m,m,\infty)$ can be found in~\cite{K}, \cite{KPTh}. 

\myskip
Part of this work was done during the stay at IHES in spring~2004.
I am grateful to IHES for its hospitality and support.
I would like to thank Martin Deraux and Richard Schwartz for useful conversations related to this work.

\section{Non-Discreteness Proof}

For fixed \ppp\ the space of complex hyperbolic triangle groups is of real dimension one.
We now describe a parameterisation of the space of complex hyperbolic triangles in~$\chp$ by means of an angular invariant $\al$. See section~3 in~\cite{P} for details.

\myskip
Let $c_k$ be the normalised polar vector of the complex geodesic $C_k$.
Let $r_k=|\<c_{k-1},c_{k+1}\>|$.
If the complex geodesics $C_{k-1}$ and $C_{k+1}$ meet at the angle $\varphi_k$, then $r_k=\cos\varphi_k$.
We define the {\it angular invariant} $\al$ of the triangle $(C_1,C_2,C_3)$ by
$$\al=\arg\left(\prod_{k=1}^3\<c_{k-1},c_{k+1}\>\right).$$
A complex hyperbolic triangle in $\chp$ is determined uniquely up to isometry
by the three angles and the angular invariant~$\al$ (compare proposition~1 in~\cite{P}). 
Let $\iota_k=\iota_{C_k}$ be the complex reflection in the complex geodesic $C_k$. 

\myskip
Let $\phi:\Gappp\to\PU(2,1)$ be a complex hyperbolic triangle group representation
and $G:=\phi(\Gappp)$ the corresponding complex hyperbolic triangle group.
Assume that $\ga$ is an element of infinite order in $\Gappp$
and that its image $\phi(\ga)$ in $G$ is regular elliptic.
Then there are two cases,
either $\phi(\ga)$ is of finite order, then $\phi$ is not injective,
or $\phi(\ga)$ is of infinite order, then $\phi$ is not discrete
because the subgroup of $G$ generated by $\phi(\ga)$ is not discrete.

\myskip
We shall show, that if the element $\iota_1\iota_2\iota_3$ is regular elliptic, then it is not of finite order, hence the corresponding triangle group is not discrete. 

\myskip
This statement was proved in~\cite{Sch2001:dented} for ideal triangle groups, \ie groups of type $(\infty,\infty,\infty)$.
The statement for $(m,m,\infty)$-triangle groups was formulated in~\cite{W} (Lemma 3.4.0.19), but the proof there had a gap.

\begin{thm}
An $(m,m,\infty)$-triangle group is not discrete if the product of the three generators is regular elliptic.
\end{thm}

\begin{proof}
We assume that the element~$\iota_1\iota_2\iota_3$ is regular elliptic of finite order.
Let~$\tau\ne-1$ be the trace of the corresponding matrix in~$\SU(2,1)$.
The eigenvalues of this matrix are then three roots of unity with product equal to~$1$.
Hence
$$\tau=\om_n^{k_1}+\om_n^{k_2}+\om_n^{k_3}$$
for some $k_1$, $k_2$, and~$k_3$ with $k_1+k_2+k_3=0$. 
Here $\om_n=\exp(2\pi i/n)$ and $n$ is taken as small as possible.
On the other hand, the trace $\tau$ can be computed (see section~8 in~\cite{P}) as
$$\tau=8r_1r_2r_3e^{i\al}-(4(r_1^2+r_2^2+r_3^2)-3).$$
Let
$$r=\cos\left(\frac{\pi}{m}\right).$$
For $(m,m,\infty)$-groups we have $r_1=r_2=r=\cos(\pi/m)$ and $r_3=1$,
hence
$$\tau=(8r^2)e^{i\al}-(8r^2+1).$$
This equation implies that the complex number $\tau\ne-1$ lies on the circle with
center in~$-(8r^2+1)$ and radius~$8r^2$, or in other words $\tau$ satisfies
the equation
$$(\tau+(8r^2+1))\cdot(\bar\tau+(8r^2+1))=|\tau+(8r^2+1)|^2=(8r^2)^2.$$
This implies in particular
$$\Re(\tau)<-1.$$
Let $N$ be the least common multiple of~$n$ and~$2m$.
Let $\si_k$ be the homomorphism of~$\q[\om_N]$ given by $\si_k(\om_N)=\om_N^k$.
For $k$ relatively prime to~$n$ the restriction of $\si_k$ to~$\q[\om_n]$ is a Galois automorphism.

\begin{lem}
\label{hom-sigma-k}
Let $\tau=\om_n^{k_1}+\om_n^{k_2}+\om_n^{k_3}$ be the trace of the matrix of~$\iota_1\iota_2\iota_3$, where~$n$ is taken as small as possible.
Then~$\si_k(\tau)$ satisfies the equation
$$|\si_k(\tau)+\si_k(8r^2)+1|=\si_k(8r^2).$$
This implies in particular
$$\Re(\si_k(\tau))\le-1.$$
\end{lem}

\begin{proof}
We have
$$\tau\in\q[\om_n]\subset\q[\om_N]$$
and
$$
  2r
  =2\cos\left(\frac{\pi}{m}\right)
  =\om_{2m}+\bar{\om}_{2m}
  \in\q[\om_{2m}]
  \subset\q[\om_N],
$$
hence the equation $|\tau+(8r^2+1)|=8r^2$ is defined in~$\q[\om_N]$.
The homomorphism~$\si_k$ commutes with complex conjugation and hence maps real numbers to real numbers.
Applying the homomorphism~$\si_k$ to the equation 
$$(\tau+(8r^2+1))\cdot(\bar\tau+(8r^2+1))=(8r^2)^2$$
we obtain
$$(\si_k(\tau)+\si_k(8r^2+1))(\si_k(\bar\tau)+\si_k(8r^2+1))=(\si_k(8r^2))^2.$$
This equation is equivalent to
$$|\si_k(\tau)+\si_k(8r^2)+1|^2=(\si_k(8r^2))^2.$$
Since $\si_k(2r)$ is a real number,
the number $\si_k(8r^2)=2(\si_k(2r))^2$ is a non-negative real number.
Hence $\si_k(\tau)$ satisfies the equation
$$|\si_k(\tau)+\si_k(8r^2)+1|=\si_k(8r^2).$$
This equation means that the complex number $\si_k(\tau)$ lies on the
circle with center in $-(\si_k(8r^2)+1)<0$ and radius $\si_k(8r^2)\ge0$.
This implies in particular
$$\Re(\si_k(\tau))\le-1.\qedhere$$
\end{proof}

\begin{lem}
\label{sph-inequality}
Let $\tau=\om_n^{k_1}+\om_n^{k_2}+\om_n^{k_3}$ be the trace of the matrix of~$\iota_1\iota_2\iota_3$, where~$n$ is taken as small as possible.
For~$i\in\{1,2,3\}$, let
$$d_i=\frac{n}{\gcd(k_i,n)},$$
where $\gcd$ is the greatest common divisor.
Then
$$\frac{1}{\varphi(d_1)}+\frac{1}{\varphi(d_2)}+\frac{1}{\varphi(d_3)}\ge1.$$
\end{lem}

\begin{proof}
According to Lemma~\ref{hom-sigma-k}, 
$$\Re(\si_k(\tau))\le-1$$
for any homomorphism~$\si_k$.
Summing over all $k\in\{1,\dots,n-1\}$ relatively prime to~$n$ we obtain
$$\Re\left(\sum\limits_{{1\le k<n}\atop (k,n)=1}\,\si_k(\tau)\right)<-\varphi(n)$$
and hence
$$\left|\sum\limits_{{1\le k<n}\atop (k,n)=1}\,\si_k(\tau)\right|>\varphi(n).$$
Here $\varphi$ is the Euler $\varphi$-function.
The root of unity $\om_n^{k_i}$ is a primitive $d_i$-th root of unity.
The sum of all $d_i$-th primitive roots of unity is in~$\{-1,0,1\}$,
and hence is bounded by~$1$.
The map $(\z/n\z)^*\to(\z/d_i\z)^*$ induced by the map $\z/n\z\to\z/d_i\z$ is surjective,
and the preimage of any element in $(\z/d_i\z)^*$ consists of
$\varphi(n)/\varphi(d_i)$ elements.
Hence we obtain the inequality 
$$\big|\sum\limits_{{1\le k<n}\atop (k,n)=1}\,\si_k(\om_n^{k_i})\big|\le\frac{1}{\varphi(d_i)}\cdot\varphi(n)$$
for $i\in\{1,2,3\}$.
From the inequalities
\begin{align*}
  \varphi(n)
  &<\big|\sum\limits_{{1\le k<n}\atop (k,n)=1}\,\si_k(\tau)\big|\\
  &=\big|\sum\limits_{{1\le k<n}\atop(k,n)=1}\,\si_k(\om_n^{k_1}+\om_n^{k_2}+\om_n^{k_3})\big|\\
  &\le\left(\frac{1}{\varphi(d_1)}+\frac{1}{\varphi(d_2)}+\frac{1}{\varphi(d_3)}\right)\cdot\varphi(n)
\end{align*}
it follows
$$\frac{1}{\varphi(d_1)}+\frac{1}{\varphi(d_2)}+\frac{1}{\varphi(d_3)}>1.$$
\end{proof}

\noindent
The inequality
$$\frac{1}{\varphi(d_1)}+\frac{1}{\varphi(d_2)}+\frac{1}{\varphi(d_3)}>1$$
implies that the triple~$(\varphi(d_1),\varphi(d_2),\varphi(d_3))$ is equal (up to permutation) to one of the triples
$$(1,?,?),\quad (2,2,?),\quad (2,3,3),\quad (2,3,4),\quad (2,3,5).$$
But for the Euler $\varphi$-function we have~$\varphi(x)=1$ for~$x\in\{1,2\}$, $\varphi(x)=2$ for $x\in\{3,4,6\}$ and $\varphi(x)\ge4$ for all other positive integers~$x$.
\begin{enumerate}[$\bullet$]
\item
The triples~$(2,3,3)$, $(2,3,4)$, $(2,3,5)$ cannot occur since~$\varphi(x)\ne3$ for any integer~$x$.
\item
Let~$\varphi(d_i)=1$ for some~$i\in\{1,2,3\}$.
Without loss of generality we can assume that~$\varphi(d_1)=1$.
Then~$d_1\in\{1,2\}$, therefore~$(k_1,n)\in\{n/2,n\}$ and~$k_1\equiv0,n/2\mod~n$.
Hence~$\om_n^{k_1}\in\{1,-1\}$.
If~$k_1\equiv0$ then~$k_2+k_3\equiv0$.
Let~$k=k_2$, then~$k_3\equiv-k$ and
$$\tau=\om_n^{k_1}+\om_n^{k_2}+\om_n^{k_3}=1+\om_n^k+\om_n^{-k}=1+2\cos(2\pi k/n)$$
and~$\Re(\tau)=1+2\cos(2\pi k/n)\ge-1$ in contradiction to~$\Re(\tau)<-1$.
If~$k_1\equiv n/2$ then~$k_2+k_3\equiv-n/2$.
Let~$k=k_2$, then~$k_3\equiv-k-n/2$ and
$$\tau=\om_n^{k_1}+\om_n^{k_2}+\om_n^{k_3}=-1+\om_n^k-\om_n^{-k}=-1+2i\sin(2\pi k/n)$$
and~$\Re(\tau)=-1$ in contradiction to~$\Re(\tau)<-1$.
\item
If~$\varphi(d_i)=\varphi(d_j)=2$ for~$i,j\in\{1,2,3\}$, $i\ne j$, then~$d_i,d_j\in\{3,4,6\}$,
therefore~$(k_i,n)\in\{n/6,n/4,n/3\}$ and~$k_i\equiv\pm n/6,\pm n/4,\pm n/3\mod~n$.
Hence~$\om_n^{k_i}\in\{\al^{\pm2},\al^{\pm3},\al^{\pm4}\}$, where~$\al=\om_{12}=\exp(2\pi i/12)$, and
$$\tau=\al^p+\al^q+\al^r,\quad p+q+r=0,\quad p,q\in\{\pm2,\pm3,\pm4\}.$$
Using~$\Re(\tau)<-1$ and~$\Re(\al^r)\ge-1$ we obtain
$$\Re(\al^p+\al^q)=\Re(\tau)-\Re(\al^r)<-1+1=0.$$
Since~$\Re(\al^{\pm2})=\frac12$, $\Re(\al^{\pm3})=0$ and~$\Re(\al^{\pm4})=-\frac12$,
we can only have~$\Re(\al^p+\al^q)<0$ if~$\al^p+\al^q=\al^{\pm3}+\al^{\pm4}$ or~$\al^p+\al^q=\al^{\pm4}+\al^{\pm4}$.
Out of these cases, we easily check that~$\Re(\al^p+\al^q+\al^{-p-q})<-1$ holds only if~$\al^p=\al^q=\al^4$ or~$\al^p=\al^q=\al^{-4}$,
i.e.\ if~$\tau=3\al^{\pm4}$, but then a suitable homomorphism~$\si_k$ has the property~$\Re(\si_k(\tau))>-1$ in contradiction to Lemma~\ref{hom-sigma-k}.
\end{enumerate}
\end{proof}


\bibliographystyle{amsalpha}
\bibliography{traces-add}

\providecommand{\bysame}{\leavevmode\hbox to3em{\hrulefill}\thinspace}
\providecommand{\MR}{\relax\ifhmode\unskip\space\fi MR }
\providecommand{\MRhref}[2]{%
  \href{http://www.ams.org/mathscinet-getitem?mr=#1}{#2}
}
\providecommand{\href}[2]{#2}
\begin{thebibliography}{KPT09}

\bibitem[Gol99]{G99}
William~M. Goldman, \emph{{Complex Hyperbolic Geometry}}, Oxford University
  Press, 1999.

\bibitem[Kam07]{K}
Shigeyasu Kamiya, \emph{{Remarks on Complex Hyperbolic Triangle Groups}},
  {Complex analysis and its applications}, OCAMI Stud., vol.~2, Osaca Munic.\
  Univ.\ Press, Osaca, 2007, pp.~219--223.

\bibitem[KPT09]{KPTh}
Shigeyasu Kamiya, John~R. Parker, and James~M. Thompson, \emph{{Non-Discrete
  Complex Hyperbolic Triangle Groups of Type $(n,n,\infty;k)$}}, 2009.

\bibitem[Par03]{P03}
John~R. Parker, \emph{{Notes on Complex Hyperbolic Geometry}}, preliminary
  version (July 11, 2003), 2003.

\bibitem[Pra05]{P}
Anna Pratoussevitch, \emph{{Traces in Complex Hyperbolic Triangle Groups}},
  Geometriae Dedicata \textbf{111} (2005), 159--185.

\bibitem[Sch01]{Sch2001:dented}
Richard~E. Schwartz, \emph{Ideal triangle groups, dented tori and numerical
  analysis}, Annals of Math. \textbf{153} (2001), 533--598.

\bibitem[WG00]{W}
Justin Wyss-Gallifent, \emph{{Complex Hyperbolic Triangle Groups}}, Ph.D.
  thesis, University of Maryland, 2000.

\end{thebibliography}

\end{document}